\documentstyle[12pt]{article}
\newcommand{\argmin}{\mathop{\rm argmin}\limits}

\newcommand{\const}{\mathop{\rm const}\limits}

\newcommand{\mod}{\mathop{\rm mod}\limits}

\newcommand{\Law}{\mathop{\rm Law}\limits}

\newcommand{\Ent}{\mathop{\rm Ent}\limits}

\begin{document}

\begin{center}

{\bf ADAPTIVE OPTIMAL REGULARIZATION OF THE LINEAR} \\

\vspace{2mm}

{\bf ILL POSED INTEGRAL EQUATIONS.} \\

\vspace{3mm}

{\bf A statistical nonparametric asymptotical approach.}\\

\vspace{3mm}

   E.Ostrovsky and L.Sirota, {\sc Israel}. \\

\vspace{2mm}

{\it Department of Mathematics and Statistics, Bar-Ilan University,
59200, Ramat Gan.}\\
e \ - \ mails: galo@list.ru; \ sirota@zahav.net.il \\

\vspace{3mm}
                     {\sc Abstract.}\\

\end{center}
\vspace{2mm}

  We construct an adaptive asymptotically  optimal in order in the $ L(2) $ sense
 a solution (estimation) of an integral linear equation  of a first kind and energy of this solution with the confidence region building, also adaptive. \par

\vspace{2mm}

 \par

\vspace{2mm}

\par

\vspace{3mm}

\section{Statement of problem.}
\vspace{3mm}

 Let us consider the following linear integral equation relative a unknown
 Riemann integrable function $ f: [0, 1] \to (-\infty,\infty) $ of a first kind
 (ill posed problem, see \cite{Tikhonov1}):

 $$
  \int_0^{1} R(t_i-s)f(s)ds + \sigma \epsilon_i = y(t_i), \eqno(1)
 $$
or equally
$$
g(t_i) + \sigma \epsilon_i = y(t_i), \ g(t) \stackrel{def}{=} R*f(t),\eqno(2)
$$
with the data $ y(t_i) $ obtained in the uniform set $ t_i = i/n, \ i = 1,2, \ldots,n $
with the random errors of measurement (noise) $ \sigma \epsilon_i. $ \par
 Here $ \{ \sigma \cdot \epsilon_i \}, \ {\bf Var}(\epsilon_i)=1 $ are centered:
 $ {\bf E} \epsilon_i = 0 $ independent random variables, which may be defined on
 some probabilistic space $ (\Omega, \cal{M}, {\bf P}) $ with an expectation $ {\bf E } $
 and variance $ {\bf Var}, $ errors of the measurements at the points $ t_i, \ \sigma = \const \in (0,\infty); \ R(\cdot) $ is the kernel, which may be generalized function; the arithmetical operation under the arguments $ t,s $ are understood modulo $ 1 $
 (periodicity). \par
  The consistent $ 1/\sqrt{n} $ estimation $ \sigma(n) $ of the variable $ \sigma $
is obtained in the article \cite{Golubev1} by means of the Residual Sum of Squares 
(RSS) method; indeed,

$$
\sigma^2(n) = \sum_{j=1}^n (y(t_i) - f_0(n, t_i))^2/(n-1),
$$
where $ f_0(n, t) $ is some {\it preliminary} estimation of the function $ f; $ for instance, $ f_0(n, t) $ may be obtained by means of the Minimum Square Estimation (MSE) method.\par
 More exact estimation of the variance of error $ \sigma^2 $  see in article    \cite{Maruyama1}; in this  work are offered and investigated the so-called supereffective Generalized Bayes Estimators for $ \sigma^2. $ \par
  We will suppose therefore for simplicity that the value $ \sigma $ is known.\par
 We consider the asymptotical statement of problem: $ n \to \infty. $ \par

\vspace{3mm}
 Notice that in the case $ R(t)= \delta(t), $ where $ \delta(t) $ denotes the Dirac
  delta-function, the problem (1) coincides with the classical {\it regression problem } of nonparametric statistics; the optimal in the $ L_2(\Omega\times [0,1]) $ sense
 adaptive algorithms for solving (estimation, on the statistical language)
 of the function $ f(\cdot) $ is described in many publication; see, for
 example, articles
\cite{Bobrov1}, \cite{Brown1}, \cite{Brown2}, \cite{Cai1}, \cite{Cai2},
\cite{Donoho1}, \cite{Donoho2}, \cite{Donoho3}, \cite{Donoho4},
\cite{Golubev1}, \cite{Juditsky1}, \cite{Lepsky1}, \cite{Ostrovsky1},
\cite{Ostrovsky2}, \cite{Ostrovsky3}, \cite{Ostrovsky4}, \cite{Pinsker1}, \cite{Polyak1}, \cite{Polyak2}; monograph \cite{Tsybakov1} and reference therein. \par
\vspace{3mm}
{\bf Our goal is to offer and investigate an adaptive asymptotically
as $ n \to \infty $ optimal in the $ L_2(\Omega\times [0,1]) $ sense estimation
of the unknown function $ f(\cdot) $ based on the observations } $ \{ y(t_i), \ i = 1,2,\ldots,n \}. $ \par

\vspace{3mm}

\section{Denotations. Assumptions.}
\vspace{3mm}
 {\bf 1.} Let  $ X = [0,1], t \in X, T = \{ \phi_k(t) \} $ be the classical
 complete normalized trigonometrical system on the set X:
 $ \phi_1(t) = 1, \phi_2(t) = \sqrt{2}\cos(2\pi t), \phi_3(t) = \sqrt{2}\sin(2\pi t), $
 $$
  \phi_4(t)= \sqrt{2} \cos(4 \pi t), \phi_5(t) = \sqrt{2} \sin(4 \pi t),
\phi_6(t)= \sqrt{2} \cos(6 \pi t), \phi_7(t) = \sqrt{2} \sin(6 \pi t),  \ldots.
 $$
   Define for the measurable function $ g:[0,1] \to R $ from the equality (2)
 $$
g(t) = \sum_{k=1}^{\infty} c(k)\phi_k(t);  \eqno(3)
$$
here
$$
 \ c(k) = \int_{0}^1 g(t) \ \phi_k(t) \ dt
$$
 be the Fourier coefficients of a function $ g(\cdot) $ over the system $ T. $ \par
 {\bf 2.} Let us write the formal Fourier expansion for the kernel $ R(\cdot): $
 $$
 R(t) = \sum_{k=1}^{\infty} \phi_k(t)/w(k),
 $$
 then formally
 $$
 f(t)= \sum_{k=1}^{\infty} c(k) w(k) \phi_k(t).
 $$
{\sc We assume }
$$
\inf_k |w(k)| > 0, \ \exists  \theta = \const \in [0,\infty), \ |w(k)| \asymp
k^{\theta}. \eqno(4)
$$
It is evident that in the case $ \theta = 0 $ the sequence $ \{  |w(k)| \} $ is
bilateral bounded:
$$
 0 <  \inf_k |w(k)| \le \sup_k |w(k)| < \infty.
$$
 Notice that under the conditions (4) and $ y(\cdot) \in L(2) $ the limit
 as $ \epsilon \to 0+ $ and $ n \to \infty $  equation, i.e.  an equation

 $$
 R*f(t) = y(t), \ t \in (0,1)
 $$
 has an unique a.e. solution $ f = f(t). $ \par

{\bf 3.} Let us denote

$$
S(N) = \sum_{k=1}^N w^2(k), \ S_2(N) = \sum_{k=N+1}^{2N} w^4(k),
$$

$$
\rho(N) = \sum_{k=N+1}^{\infty} c^2(k) w^2(k), \ \rho_2(N) =
\sum_{k=N+1}^{2N} c^2(k) w^4(k),
$$

$$
A(N,n) = \sigma^2 S(N)/n + \rho(N), \ N^+ = N^+(n) =
$$

$$
 \min \left[ \Ent(n/\log(n+8)), \Ent((n/\log(n+8))^{1/(2 \theta)}) \right],
$$
$ \Ent[z] $ denotes here the integer part of the real positive variable $ z; $
$$
A^*(n) = \min_{N \in [1,N^+]} A(N,n), \ N_0 = \argmin_{N \in [1,N^+]} A(N,n).
$$
 Note that if $ f \in B(w), $ then
 $$
 \lim_{N \to \infty} \rho(N) = 0.
 $$
 Therefore
 $$
 \lim_{n \to \infty} A^*(n) = 0,
 $$
as long as
$$
A^*(n) \le \sigma^2 \ n^{-1} \ S \left(\Ent(n^{1/(4\theta)}) \right) +
\rho\left(\Ent(n^{1/(4\theta)}) \right) \to 0, \ n \to \infty.
$$

{\sc Assumptions:}

$$
0 <  \gamma_- \stackrel{def}{=} \underline{\lim}_{N \to \infty}\rho(2N)/\rho(N) \le
\overline{\lim}_{N \to \infty}\rho(2N)/\rho(N)  \stackrel{def}{=} \gamma_+  < 1, \eqno(5a)
$$
{\sc and we suppose for the construction of confidence region for the function
$ f $ that}

$$
\exists \ \lim_{N \to \infty}\rho(2N)/\rho(N) \stackrel{def}{=} \gamma \in (0,1).
\eqno(5b).
$$

 {\sc Another conditions:}
$$
 \Gamma^+ \stackrel{def}{=} \overline{\lim}_{N \to \infty} S(2N)/S(N)
 \in (1, \infty), \eqno(6a)
$$
{\sc or more strictly}
$$
\exists \ \Gamma \stackrel{def}{=}
{\lim}_{N \to \infty} S(2N)/S(N), \  \Gamma \in (1, \infty); \eqno(6b)
$$
and we denote

$$
U = U(\gamma_-, \Gamma_-) = \min \left[ (1-\gamma_-), \ (\Gamma_- - 1) \right].
$$
 The condition (6a) is satisfied if $ |w(k)| \asymp k^{\theta}, \ k \to \infty;  $
the condition (6b) is satisfied if $ |w(k)| \sim k^{\theta}, \ k \to \infty. $ \par

 This conditions are direct analogues of the notorious $\Delta_2 $ condition in the
theory of Orlicz's spaces.\par
The conditions (5a) and (6a) are satisfied, for instance, when as $ k \to \infty  $
$$
|c(k)| \asymp \ k^{-\Delta}, \ \Delta = \const, \ \Delta > 2 \theta + 1/2. \eqno(7)
$$
 It may be considered analogously a more general case  when

$$
\inf_k |w(k)| > 0, \ \exists  \theta = \const \in [0,\infty), \ |w(k)| \asymp
L_1(k) \ k^{\theta},
$$
$$
|c(k)| \asymp L_2(k) \ k^{-\Delta}, \ \Delta = \const, \ \Delta > 2 \theta + 1/2,
$$
where $ L_1(k), \ L_2(k) $ are slowly varying as $ k \to \infty $ functions. The detail
investigation of such a functions see in a book \cite{Seneta1}.\par
{\bf 4.} We denote by  $ c(k,n) $ usually consistent $ 1/\sqrt{n} $  estimations of a Fourier coefficients $ c(k) $ of the function $ g(\cdot) $
 based on the sample  $ y(t_i), i = 1,2,\ldots,n, $ namely:
 $$
 c(k,n) = n^{-1}\sum_{i=1}^n y(t_i) \cdot \phi_k(t_i), \eqno(8)
 $$
and define correspondingly

$$
\tau(N,n) = \sum_{k=N+1}^{2N}w^2(k) c^2(k,n), \ \tau^*(n) = \min_{N \in [1,N^+]} \tau(N,n), \eqno(9)
$$

$$
M(n) = \argmin_{N \in [1,N^+]}\tau(N,n). \eqno(10)
$$
Note that the variables $ \tau(N,n), \tau^*(n) $ and $ M(n) $ are random variables
which dependent on the source data $ \{ y(t_i)\}.$

The consistent as $ n \to \infty $ estimation $ \gamma(n) $ of the parameter
$ \gamma $ under condition (5b) is described in \cite{Ostrovsky4}, chapter 5, section
13. Namely, we define
$$
G = G(n) \stackrel{def}{=} \Ent(\exp(\sqrt{\log n}));
$$
then the value $ \gamma $ may be consistent under condition (5b)
estimated by means of statistic $ \gamma(n) $ as follows:
$$
\gamma(n)= \frac{\tau(4G) - 2 \tau(2G)}{\tau(2G) - 2 \tau(G)}.
$$

{\bf 5.} {\sc Further, we suppose  for the construction of asymptotical confidence
region that}

$$
\sup_i {\bf E} (\epsilon_i)^4  <  \infty, \eqno(11a)
$$
{\sc and we suppose for the construction of non-asymptotical confidence
region that} $ \exists q,Q = \const > 0,  $
 $$
\sup_i \max \left[ {\bf P}(\epsilon_i > u),{\bf P}(\epsilon_i < -u) \right] \le
  \exp \left(-  (u/Q)^q \right), \ u > 0. \eqno(11b)
  $$
 We will write, as ordinary, in some concrete passing to the limit, for instance,
 $ n \to \infty, $
$$
X(n) \sim Y(n)  \Leftrightarrow \lim_{n \to \infty} X(n)/Y(n) = 1,
$$

$$
X(n) \asymp Y(n)  \Leftrightarrow \inf_{n} X(n)/Y(n) \le \sup_{n} X(n)/Y(n)
< \infty.
$$
all the relations between the random variables, for example, passing to the limit,
are understood with probability one $ \ (\mod {\bf P}). $ \par

\vspace{3mm}

\section{Modular spaces.}
\vspace{3mm}

 Let again  $ X = [0,1], x \in X, T = \{ \phi_k(x) \} $ be the classical complete
 normalized trigonometrical system.  Define for the measurable function
 $ g:[0,1] \to R $
 $$
g(x) = \sum_{k=1}^{\infty} c(k)\phi_k(x);
$$
here
$$
 \ c(k) = \int_{0}^1 f(x) \ \phi_k(x) \ dx
$$
 be the Fourier coefficients of a function $ f(\cdot) $ over the system $ T. $ \par
 Let also $ w = \{ w(k) \} $ be positive number sequence (weight). We introduce the following modular space
$ B(p,w) = B(T, p,w), $ which will called {\it modular weight space,} consisting on
all the (measurable) functions $ \{g \} $ with finite norm

 $$
||g||B(p,w) = ||f||_{p,w} \stackrel{def}{=}
\left[ \sum_{k=1}^{\infty} w^p(k)|c(k)|^p  \right]^{1/p}, \ p \in (1,\infty). \eqno(12)
$$
For instance, the Sobolev's spaces $ W_2^m = W_2^m[0,1], m = 0,1,2,\ldots $
consisting on all the periodical $ (\mod 1)  \ m-$ times differentiable functions $ f $ with finite norm

$$
||f||W_2^m = \left[(||f||L_2[0,1])^2 + (||f^{(m)}||L_2[0,1])^2 \right]^{1/2}
$$
is weighted modular space relative the system $ T. $ \par
 Of course, the $ B(w) $ space with $ w(k) = 1 $ coincides with the ordinary space
 $ L(2) $ on the set $ [0,1].$  We will write in this case

$$
||f|| = ||f||L(2) = ||f||B(1).
$$

 \vspace{3mm}

  The complete investigation of these spaces see in the book \cite{Musielak1}.\par
In the case when $ p=2 $ we will write for brevity $ B(2,w) = B(w). $ Of course, the
space $ B(w) = B(2,w) $ is (separable) Hilbert space. \par
The notion  $ "L_2(\Omega\times [0,1])" $ sense  means by definition that we consider
the following {\it loss function:}

$$
V(h(n,\cdot), f) \stackrel{def}{=}  {\bf E} ||h(n,\cdot) - f(\cdot)||^2. \eqno(13)
$$
where $h(n,\cdot) $ is arbitrary estimation of a function $ f(\cdot) $ based on the
observations $ \{ y_i \}, \ i = 1,2,\ldots,n. $\par

 {\sc We suppose in addition that }
 $$
 f \in L(2)  \Leftrightarrow ||f||^2 =
 \sum_{k=1}^{\infty} c^2(k) w^2(k) < \infty. \eqno(14a)
 $$
 The condition (14a) will be used for non-adaptive estimation of a solution $ f. $ For the construction of the adaptive estimation we need to assume the following condition (14b):
$$
f \in B(w) \Leftrightarrow ||f||^2(B,w) = \sum_{k=1}^{\infty} c^2(k) w^4(k) < \infty. \eqno(14b)
$$

\vspace{3mm}

\section{ Main result. Construction of our solution (estimation).}

\vspace{3mm}
 Let us consider a projection, or Tchentsov's estimation of a function $ g(\cdot) $
of  a view

$$
g(N,n,t) = \sum_{k=1}^N c(k,n) \phi_k(t), \ N = \const \in [1, n/3],
$$
and we construct correspondingly  the following projection estimation of a function
$ f(\cdot):$

$$
f(N,n,t) = \sum_{k=1}^N c(k,n) \ w(k) \ \phi_k(t), \ N = \const \in [1, n/3]. \eqno(15)
$$
 We find by the direct calculation using the condition $ f(\cdot) \in B(w) $  as
 $ n \to \infty: $

$$
V(f(N,n,\cdot),f(\cdot)) \sim A(N,n). \eqno(16)
$$
 If we choose
 $$
 N = N_0 \stackrel{def}{=} \argmin_{N \in [1,N^+]} A(N,n),
 $$
 then we conclude that the optimal speed of convergence for
 {\it non-adaptive estimation} $ f(N_0,n,x) $  is asymptotically $ A^*(n): $ as
$ n \to \infty $
$$
{\bf E}||f(n,N_0,\cdot)-f(\cdot)||^2 \sim A^*(n).
$$
 At the same result up to multiplicative constant
 is true for our {\it adaptive estimation,} which we built by the
 following way. Namely, let us define
 $$
\hat{f}(n,t) \stackrel{def}{=} f(n,M(n),t). \eqno(17)
 $$

 The {\it asymptotically  exact adaptive estimation} of the solution
$ f(\cdot) $ may be obtained under the conditions (5b),(6b),(11a), and (14b)
by using the so-called {\it penalty function method.} Indeed, if we define
$$
M_1(n) \stackrel{def}{=} \argmin_{N \in [1,N^+]} \tau_1(N,n), \eqno(18)
$$
$$
\tau_1(N,n) \stackrel{def}{=}  \tau(N,n) + (2 - \gamma(n)- \Gamma)
 \cdot \sigma^2(n) \cdot S(N)/n =
$$

$$
\sum_{k=N+1}^{2N} c^2(k,n) w^2(k) + (2 - \gamma(n)- \Gamma)\cdot \sigma^2(n) \cdot
S(N)/n, \eqno(19)
$$
$$
\tau^*_1(n) := \min_{N \in [1,N^+]}\tau_1(N,n). \eqno(20)
$$
 The function $ (2 - \gamma(n)- \Gamma) \cdot \sigma^2(n) \cdot S(N)/n $ is
 said to be {\it penalty function.} \par
 We introduce a new estimation $ \tilde{f}(n,t) $ of the function $ f(t) $ as follows:
$$
\tilde{f}(n,t) \stackrel{def}{=} f(n,M_1(n),t). \eqno(21)
$$
We assert that under the conditions of (5b),(6b), (11b) and (14b)

$$
 V(\tilde{f}(n,\cdot), f(\cdot)) \le A^*(n)(1+\nu(n)), \eqno(22)
$$
where $ \lim_{n \to \infty} \nu(n) = 0. $  \par

\vspace{3mm}
 {\bf Theorem 1.} We suppose that all our conditions, in particular, the conditions
 (5a), (6a), (11a), and (14a) are satisfied. Then

$$
\overline{\lim}_{n \to \infty} V(\hat{f}(n,\cdot), f(\cdot))/A^*(n) \le 1, \eqno(23a)
$$

$$
\overline{\lim}_{n \to \infty} || \hat{f}(n,\cdot) - f(\cdot)||^2/\tau^*(n) \le 1/
 U(\gamma_-, \Gamma_-) \ (\mod {\bf P}). \eqno(23b)
$$

{\bf Remark 1.} Note that the estimation $ \hat{f}(n,\cdot) $ is optimal in general case, namely, in the case when $ U(\gamma_-, \Gamma_-) \to 1-. $  See \cite{Ibragimov1}. \par

\vspace{3mm}

{\bf Theorem 2.} We have  under the conditions of (5b),(6b), (11a) and (14b)

$$
 V(\tilde{f}(n,\cdot), f(\cdot)) \le A^*(n)(1+\nu_1(n)),
 \eqno(24a),
$$

$$
\overline{\lim}_{n \to \infty} ||\tilde{f}(n,\cdot) - f(\cdot)||^2
/ \tau_1^*(n) \le 1. \eqno(24b),
$$
where $ \lim_{n \to \infty} \nu_1(n) = 0. $ \par
The proof is at the same as the proof of theorem 1, see further.\par

\vspace{3mm}

\section{Energy estimation.}

\vspace{3mm}
Let us define the {\it energy functional}, or briefly {\it energy}
 $ H = H(f) $ of the signal $ f(\cdot) $ as usually

$$
H = H(f) \stackrel{def}{=} ||f||^2(L(2)) = \int_0^1 f^2(t) dt. \eqno(25)
$$
We estimate in this section the energy functional $ H = H(f); $ we will prove that under natural conditions, in particular, $  f \in B(w), $ there exists an estimation $ H(n,f) $ which convergent to the $ H(f) $ with the speed $ 1/\sqrt{n} $
(optimality). \par
 The necessity of the condition $ f \in B(w) $ for the possibility of optimal in order estimation $ H(f) $ is proved in \cite{Efromovich1}, \cite{Klemela1}.\par
  We describe now our energy estimation $ H(n,f),$ which is based on the our adaptive estimation $ \tilde{f}(n,\cdot) $ and is different on all others such estimations, see
 \cite{Donoho4}, \cite{Efromovich1},   \cite{Klemela1} etc. \par
 Our estimation is some slight generalization of offered therein. \par
  We introduce the following energy estimation functional:

 $$
 H(n,f) = \sum_{k=1}^M c^2(k,n) w^2(k) - \sigma^2(n) S(M)/n. \eqno(26)
 $$
 We will called the function $ \sigma^2(n) S(M)/n $ as {\it anti-penalty function}.\par
{\bf Theorem 3.} Let $ f \in B(w). $ Assume in addition to the conditions
(5b), (6b). (11a)

$$
\lim_{n \to \infty} S_2(N_0(n))/n = 0, \
\lim_{n \to \infty}\rho(N_0(n))/\sqrt{n} = 0. \eqno(27)
$$
Then the statistics  $ H(n,f) $ is $ 1/\sqrt{n} $ consistent and optimal up to multiplicative constant estimation of the energy value $ H:$

$$
\lim_{n \to \infty}  n \cdot {\bf E} (H(n,f) - H(f))^2 = 4 \cdot \sigma^2 \cdot
||f||^2(B(w)). \eqno(28)
$$
Notice that it follows from the proposition (28) the optimality of energy
estimation $ H(n,f), $ see \cite{Efromovich1}.\newline
\vspace{3mm}

\section{Proofs.}

\vspace{3mm}
{\bf A. Investigation of the function estimations; theorems 1 and 2.} \par
{\bf 1.} The proof is at the same as in \cite{Ostrovsky4}, chapter 5, section 13;
see also \cite{Bobrov1}, \cite{Ostrovsky1} etc.,  where is considered the case
 $ R(t) = \delta(t), $ i.e. the case of the classical regression
problem in the nonparametric statistics. Namely, it is proved therein that as
$ n \to \infty $ the empirical Fourier coefficients $ c(k,n) $ are
common asymptotically independent and have normal  distribution with the parameters

$$
\Law(c(k,n)) = N(c(k), \sigma^2/n).
$$
 The detail proof of this assertion see, e.g. in \cite{Brown1}, \cite{Cai2}, \cite{Donoho4},  \cite{Efromovich1}, \cite{Ostrovsky1}, \cite{Tsybakov1}. \par

Therefore, we can write the following representation

$$
c(k,n) \stackrel{\sim}{=} c(k) + \sigma \zeta(k)/\sqrt{n}, \eqno(29)
$$
where the random variables $ \{ \zeta(k) \} $ common independent and have standard
normal  distribution.

{\bf 2.} We find by the direct calculations using the representation (29) for the projection estimations $ f(N,n,t) $ as $  N,n \to \infty, \ N \in [1,N^+]: $

$$
||f(N,n,\cdot) -f||^2 \sim \sigma^2 \ n^{-1} \sum_{k=1}^N \zeta^2(k) w^2(k) + \rho(N)=
$$

$$
\rho(N) + \sigma^2 \ n^{-1} \sum_{k=1}^N w^2(k) +
\sigma^2 \ n^{-1} \sum_{k=1}^N w^2(k) \ (\zeta^2(k)-1) \stackrel{def}{=}
$$

$$
A(N,n) + \beta(N,n). \eqno(30)
$$
It follows from the Law of Iterated Logarithm (LIL) that

$$
\overline{\lim}_{n \to \infty} \min_{N \in [1,N^+]} \beta(N,n)/A^*(n) = 0,
$$

$$
A(N,n)= {\bf E} ||f(N,n,\cdot) -f||^2 \sim  \rho(N) + \sigma^2 n^{-1} \ S(N),
$$
therefore
$$
\overline{\lim}_{n \to \infty} V(\hat{f}(n,\cdot), f(\cdot))/A^*(n) \le
1/(1- \gamma_+).
$$
{\bf 3.} As long as for $ n \to \infty $

$$
{\bf E} \tau(n) \ge \rho(N) (1-\gamma_-) + \sigma^2 \
n^{-1}(\Gamma_- -1)S(N) \ge U(\gamma_-, \Gamma_-)\cdot A(N,n),
$$

$$
A^*(n) \le {\bf E}\tau^*(n)/U(\gamma_-, \Gamma_-),
$$
we conclude

$$
\overline{\lim}_{n \to \infty}A^*(n)/\tau^*(n) \le 1/U(\gamma_-, \Gamma_-).
$$
This completes the proof of theorem 1. Theorem 2 is proved analogously.\par

{\bf 4.} Further, we conclude analogously by means of the condition (5a):
$$
{\bf E}[\tau(N,n)] \asymp  \rho(N) + n^{-1} \ \sigma^2 \ S(N) = A(N,n). \eqno(31)
$$

 If we take  in (31) the value $ N = N_0, $ we conclude

 $$
{\bf E} ||f(N_0,n,\cdot) -f||^2  \asymp A^*(n).
 $$
 Hence the estimation $ f(N_0,n,\cdot) $ is consistent estimation of a function
 $ f(\cdot) $ in the $ L(2) $ sense. \par
  It follows from the theorems Tchentsov \cite{Tchentsov1} and Ibragimov-Chasminsky
  \cite{Ibragimov1} that the speed of convergence $ \sqrt{A^*(n)} $  as $ n \to \infty $ is optimal in general case. But this method of
 solution need the value $ \rho(N) $ or at least its order as $ N \to \infty.$
 This solution $ f(N_0,n,\cdot) $ is called {\it non-adaptive}.\par
  The solution (estimation) which does not use any apriory information about estimating  function $ f(\cdot) $ will be called {\it adaptive.} \par
   \vspace{2mm}

{\bf 5.} We compute as $ N,n \to \infty, N \in [1, N^+] $
$$
{\bf Var}[\tau(N,n)] \asymp n^{-1} \ [S_2(N) + n^{-1} \ \rho_2(N)].
$$

 We conclude by virtue of adaptive conditions that for some  positive constant
 $ \beta > 0 $

 $$
{\bf Var}[\tau(N,n)]  \le C \ n^{-\beta} \ \cdot A(N,n).
 $$
 Further considerations are alike to the \cite{Ostrovsky4}, chapter 5, section 13;
 see also \cite{Bobrov1}, \cite{Ostrovsky1}. \par

 Note that it is proved there that under conditions (5b), (6b) that
with probability one and in the $ L_2(\Omega) $ sense
$$
\lim_{n \to \infty} M_1(n)/N_0 = 1, \
\lim_{n \to \infty}\rho(M_1(n))/\rho(N_0) = 1. \eqno(32)
$$
As a consequence: the estimation $ \tilde{f}(n,\cdot) $ is consistent estimation
of the solution $ f $ in the $ L(2) $ sense with probability one:

$$
||\tilde{f}(n,\cdot) - f(\cdot)||^2 \to 0, n \to \infty \ (\mod {\bf P}).
$$
 At the same result is true for the estimation $ \hat{f}(n,\cdot).$ \par

\vspace{3mm}

{\bf B. Investigation of the energy estimations; theorem 3.} \par
{\bf 1.}  We have for the energy estimation $ H(n,f), \ (M = M(n)): $

$$
H(n,f) \sim \sum_{k=1}^M w^2(k)[ c(k) + \sigma \zeta(k)/\sqrt{n} ]^2 - \sigma^2 S(M)/n \sim
$$

$$
\sum_{k=1}^M c^2(k) w^2(k) + 2 \sigma n^{-1/2} \sum_{k=1}^M w^2(k) c(k) \zeta(k) +
\sigma^2 n^{-1} \sum_{k=1}^M w^2(k)(\zeta^2(k) - 1) \sim
$$

$$
\sum_{k=1}^{N_0} c^2(k) w^2(k) + 2 \sigma n^{-1/2} \sum_{k=1}^{N_0} w^2(k) c(k) \zeta(k)+ \sigma^2 n^{-1} \sum_{k=1}^{N_0} w^2(k)(\zeta^2(k) - 1) =
$$

$$
H(f) - \rho(N_0) + 2\sigma \eta_1(n) + \sigma^2 \eta_2(n),\eqno(33)
$$
where
$$
\eta_1(n) = n^{-1/2} \sum_{k=1}^{N_0} w^2(k) c(k) \zeta(k),
$$

$$
\eta_2(n) = n^{-1} \sum_{k=1}^M w^2(k)(\zeta^2(k) - 1).
$$
 It follows from the condition (27) that

$$
\lim_{n \to \infty}  n \cdot {\bf E} (H(n,f) - H(f))^2 = \lim_{n \to \infty}
{\bf Var} [\sum_{k=1}^{N_0} w^2(k) c(k) \zeta(k)]=
$$

$$
{\bf Var} [\sum_{k=1}^{\infty} w^2(k) c(k) \zeta(k)]= ||f||^2(B(w)) < \infty,
$$
as long as $ ||f||(B(w)) < \infty. $ \par
 This completes the proof of theorem 3. \par

\vspace{3mm}

\section{Confidence regions.}

\vspace{3mm}
{\bf A.} Confidence interval (adaptive) for the function $ f $ in the $ L(2) $ norm. \par
For the rough building of the confidence domain, also
 adaptive, in the $ B(w) $ norm we proved the following result. \par

{\bf Theorem 4.} We assert under at the same conditions as in theorem 1

$$
\overline{\lim}_{n \to \infty} \tau^*(n)/A^*(n)\le 1/U(\gamma_-, \Gamma_-) \eqno(34)
$$
and correspondingly
$$
\overline{\lim}_{n \to \infty}|| \hat{f}(n,\cdot) - f(\cdot)||^2/ \tau^*(n)
\le 1/U(\gamma_-, \Gamma_-). \eqno(35)
$$

Therefore, we conclude: with probability tending to one as $ n \to \infty $

$$
|| \hat{f}(n,\cdot) - f(\cdot)||^2 \le 1.05 \cdot  \tau^*(n)/(1-\gamma_+)^2. \eqno(36)
$$
But in general case the value $ \gamma_+ $ is unknown; in order to construct
{\it estimable} confidence region, we need to suppose more strictly conditions
(5b), (6b), (11b) and (14b). \par

\vspace{3mm}

{\bf Theorem 5.} We assert under  the  conditions (5b), (6b), (11a) and (14b)
that the variables
$$
\Delta^2(n) \stackrel{def}{=} ||\tilde{f}(n,\cdot) - f(\cdot)||^2
$$
may be represented as follows:
$$
\Delta^2(n) = \tau^*(n) + \sigma(n) \ \sqrt{2M_1(n)/n} \times \xi(n), \eqno(37)
$$
where the sequence of a random variables $ \xi(n) $ has asymptotically standard
normal (Gaussian) $ N(0,1) $ distribution:

$$
\lim_{n \to \infty} \Law(\xi(n)) = N(0,1).\eqno(38a)
$$
The {\bf Proof} follows from the decomposition (30) and the Central Limit Theorem
(CLT) for the sum

$$
\xi(n) = [\sigma(n) \ \sqrt{2M_1(n)/n}]^{-1} \sum_{k=1}^{M_1(n)}(\zeta^2(k)-1).
$$

Notice that

$$
\lim_{n \to \infty} [\sqrt{M_1(n)/n}]/A^*(n)= 0 \ (\mod {\bf P}).
$$
 Non-asymptotical confidence interval for $ f(\cdot) $ may be constructed as in
\cite{Bobrov1} in the case  classical regression problem $ R(t) = \delta(t) $
under the condition (11b). Namely, let us denote

$$
r = r(q) = \min(q/2,2);
$$
then

$$
\sup_{n \ge 16 } \max [ {\bf P}(\xi(n)> Qu), {\bf P}(\xi(n)< -Qu)] \le
\exp \left(-C u^r \right).  \eqno(38b)
$$

\vspace{3mm}
{\bf B.} Adaptive confidence interval for the energy. \par
\vspace{3mm}
The next result may be proved analogously, by means of decomposition (32)
instead  (30) used by the proof of theorem 5. \par
 {\bf Theorem 6}. We assert under the conditions (5b), (6b), (11a), (14b) and the condition $ f \in B(w) $ that the estimation $ V(n) $ has asymptotically normal
 distribution with parameters

$$
\Law(H(n,f)) \sim N(H(f), 4\sigma^2 ||f||^2B(w)/n). \eqno(39)
$$
{\bf Remark 2.} The value $||g||^2(B(w^2) = ||f||^2B(w) $  may be estimated alike the value $ V = ||f||^2. $ \newline
{\bf Remark 3.} In the case when $ w(k) = 1, $ i.e.when $ g(t) = f(t) $ or equally
 $ R(t) = \delta(t), $ we have $ ||g||^2 B(w) = ||f||^2 = H(f); $ hence

$$
\Law(H(n,f)) \sim N(H(f), 4\sigma^2 H(f)/n).\eqno(40)
$$
 Using the Fisher's transform,  we conclude that the variable

 $$
\sqrt{n}(\sqrt{H(n,f)} -  \sqrt{H(f)})/\sigma(n) \eqno(41)
 $$
has asymptotically Gaussian standard distribution. \par

\vspace{3mm}

 \section{Concluding remarks.}

\vspace{3mm}

{\bf 1.}  The optimal consistent estimation  $ \hat{g}(n,t) $ in the $ L(2) $ sense of a function $ g(\cdot) $ ("Regression problem") offered in \cite{Bobrov1}, \cite{Ostrovsky1}, \cite{Ostrovsky4} has a view:

$$
\hat{g}(n,t) = \sum_{k=1}^{N_1(n)} c(k,n) \ \phi_k(t),
$$
where

$$
N_1(n) = \argmin_{N \in [1, \Ent(0.5 n)]} \sum_{k=N+1}^{2N} c^2(k,n).
$$
 But the estimation of a function $ f(\cdot) $ of a view

$$
\hat{f}_1(n,t) = \sum_{k=1}^{N_1(n)} c(k,n) \ w(k) \phi_k(t),
$$
based on the estimation  $ \hat{g}(n,t), $ is not optimal when
$ \lim_{k \to \infty} |w(k)| = \infty. $ \par

\vspace{3mm}

{\bf 2.} Adaptive estimation in general modular spaces.\par
In order to construct an adaptive optimal in order as $ n \to \infty $
estimation of the function $ f $ in the $ B(p,w), \ p \in (1,\infty) $
norm, we introduce the following function:

$$
 f^{(p)}(n,x):= \sum_{k=1}^{M^{(p)}(n)} w(k) c(k,n) \phi_k(x),
$$
 where
$$
M^{(p)}(n):= \argmin_{N \in \left[1,(N^+)^{1/p} \right]}\sum_{k=N+1}^{2N}
  w^p(k)|c(k,n)|^p.
$$
 More detail investigation will be publish in an another article.\par

{\bf 3.} Multidimensional case. \par
We consider in this subsection  the following multidimensional generalization
of our problem. Let  $ Z(n),  n = 16,17, \ldots  $ be a sequence of a vector-valued  sets  (plans of experiences) in the cube $ [0,1]^d, \ d = 2,3, \ldots: $

$$
Z(n) = \{ x_i = \vec{t}_i = \vec{t}_i(n), \}, \
\vec{t}_i \in [0, 1]^d.
$$

 At the points $ \vec{t}_i $ we observe the unknown signal (process, field)
$ f = f(t), \ t \in [0,1]^d $ on the background noise:
$$
y(t_i) = R*f(\vec{t}_i) + \sigma \ \epsilon_i,
$$
where the sequence noise $ \{\sigma \epsilon_i \}, \ $ is the sequence of errors of measurements, is the sequence of independent (or weakly dependent) centered:
 $ {\bf E } \xi_i = 0 $ normalized: $ {\bf Var }(\epsilon_i) = 1 $ random variables, $ \sigma = \const > 0 $ is a standard deviation of errors. \par
 The investigation of this problem in the case $ p = 2 $ and $ R(t)= \delta(t) $
 see in \cite{Ostrovsky5}; we will only emphasis here the neediness the using optimal planing of experience, in other words, experience design. \par
{\bf 4.} Example.\par
Let again
$$
|c(k)| \sim C_1 k^{-\Delta}, \ |w(k)| \sim C_2 k^{\theta}, \ \Delta, \theta =
\const > 0.
$$
 If

 $$
 \Delta > 2 \theta + 1,
 $$
then as $ n \to \infty $

$$
V(\tilde{f}(n,\cdot),f) \sim C_3(\Delta,\theta)
n^{-(2\Delta - 2 \theta - 1)/(2 \Delta)}
$$
and
$$
{\bf E}|V(n,f) - V(f)|^2 \sim C_4(\Delta,\theta) n^{-1}.
$$


\vspace{3mm}
\begin{thebibliography}{99}

\vspace{3mm}

\bibitem{Bobrov1}
Bobrov P.B., Ostrovsky E.I. {\it Confidence intervals by adaptive regression
estimation.} In: Probability and Statistics, Sant Petersburg, (1997), {\bf 2}, B. 244, p. 28-45; in Russian.
\bibitem{Brown1}
Brown L.D. and Low M.G.{\it Asymptotic equivalence of nonparametric regression and white noise. } Annals Stat., {\bf 24,} (1996), p. 2384-2398; MR 14259518.
\bibitem{Brown2}
Brown L.D., Cai T.T. and Zhou H.H. {\it Robust nonparametric estimation via wavelet
median regression.}
Annals Stat., {\bf 36,} (2008), p. 2055-2084; MR 15654307.
\bibitem{Cai1}
Cai T.T. {\it Adaptive wavelet estimation. A block thresholding and oracle inequality approach.}
Annals Stat., {\bf 27,} (1999), p. 898-924; MR 1724038.
\bibitem{Cai2}
Cai T.T. and Harrison H. {\it Asymptotic equivalence and adaptive estimation for robust
nonparametric  regression.}
Annals Stat., {\bf 37,} N. 6, Dezember,(2009), p. 3204-3236.
\bibitem{Donoho1}
 Donoho D.L.  {\it Wedgelets: nearly minimax estimation of edges.}  Annals of  of Statist., 1999, v. 27 b. 3 pp. 859 - 897.
\bibitem{Donoho2}
 Donoho D.L.  {\it Unconditional bases are optimal bases for data compression and for statistical estimation.} Applied Comput. Harmon. Anal., 1996, v. 3 pp. 100 - 115.
\bibitem{Donoho3}
Donoho D.L., Johnstone I.M., Kerkyacharian G., and Picard D. {\it Wavelet shrinkage, asymptotic (with discussion).} Journal of Royal Statistic; Soc, Ser. B(57), (1995), p. 301-369, MR 1323344.
\bibitem{Donoho4}
Donoho D.L. and Nussbaum M.{\it Minimax quadratic estimation of a quadratic
functional.} J. Complexity, (1990), {\bf 6}, p. 290-323, MR 1081043.
\bibitem{Efromovich1}
Efromovich S.Y. and Low M. {\it On the estimation of quadratic functionals.}
Annals Stat., {\bf 24,} (1996), p. 1106-1125; MR 1401840.
\bibitem{Golubev1}
 Golubev G., Nussbaum M.  {\it Adaptive spline Estimations in the nonparametric Regression Model.}  Theory Probab. Appl., 1992, v. 37 $ N^o $ 4, 521 - 529, in Russian.
\bibitem{Ibragimov1}
Ibragimov I.A., Chasminsky R.Z. {\it On the boundaries of quantity of nonparametrical
regression estimation.} Theory Probab. Appl., (1982), V. 21 No 1, p. 81-94, in Russian.
\bibitem{Juditsky1}
 Juditsky A., Nemirovsky A. {\it Nonparametric denoising Signals of Unknown
Local Structure, II: Nonparametric Regression Estimation.} Electronic Publication, arXiv:0903.0913v1 [math.ST] 5Mar2009.
\bibitem{Klemela1}
Klemela J. {\it Sharp adaptive estimations of quadratic functionals.}
Probab. Theory Related Fields, {\bf 134}, (2006), p. 539-564, MR 2214904.
\bibitem{Lepsky1}
Lepsky O.V. {\it On a problem of adaptive estimation in white Gaussian noise.}
Theory Probab. Appl., {\bf 35}, (1990), p. 454-466, MR 1091202, in Russian.
\bibitem{Maruyama1}
Maruyama Y., W.E. Strawderman. {\it  Improved robust Bayes estimators of the 
error variance in linear models.}
Electronic Publications, arXiv:1004.023v1 [math.ST] 1 Apr 2010.
\bibitem{Musielak1}
E.L.Musielak. {\it Modular spaces.} Kluvner Verlag, (2003), Dorderecht.
\bibitem{Ostrovsky1}
 Ostrovsky E., Sirota L. {\it Universal adaptive estimations and confidence  intervals in the non-parametrical statistics.}
 Electronic Publications, arXiv.mathPR/0406535 v1 25 Jun 2004.
\bibitem{Ostrovsky2}
 Ostrovsky E., Zelikov Yu. {\it Adaptive Optimal Nonparametric Regression and Density Estimation based on Fourier - Legendre Expansion. }
 Electronic Publication, arXiv:0706.0881v1 [math.ST] 6 Jun 2007.
\bibitem{Ostrovsky3}
Ostrovsky E.,Sirota L. {\it Optimal adaptive nonparametric denoising of
multidimensional-time signal.} Electronic Publication,
arXiv:0809.30211v1 [physics.data-an] 17 Sep 2008.
\bibitem{Ostrovsky4}
 Ostrovsky E.I. {\it Exponential Estimations for Random Fields.}
Moscow - Obninsk, OINPE, (1999), in Russian.
\bibitem{Ostrovsky5}
OstrovskyE., Rogovee E. and Sirota L. {\it Optimal Adaptive Signal Detection
and Measurement.} In: Abstracts of the International Symposium on {\sc stochastic  models in reliability engineering, life sciences and operations management,}
Beer Sheva, Israel, (2010), p. 175.
\bibitem{Pinsker1}
Pinsker M.S. {\it Optimal filtering of square integrable signals in Gaussian white
noise.} Problem Inform. Transmission, {\bf 27}, (1980), p. 120-133, in Russian.
\bibitem{Polyak1}
 Polyak B., Tsybakov A. $ C_p - $ {\it criterion in projective Estimation of
Regression.} Theory Probab. Appl., (1990), v. 35 $ N^o $ 2, 293-306, in Russian.
\bibitem{Polyak2}
Polyak B., Tsybakov A. {\it A family of asymptotically optimal Methods for selection
the order of projective Estimation of Regression.}
Theory Probab. Appl., (1992), v. 37 $ N^o $ 3, 471-485, in Russian.
\bibitem{Seneta1}
Seneta E. {\it Regularly Varying Functions.} Lecture Notes in Mathematics, Springer Verlag, (1976), Berlin-Heidelberg-New York.
\bibitem{Tchentsov1}
Tchentsov N.N. {\it Statistical resolving rules and optimal conclusions.} Moscow,
Nauka, (1972),  in Russian.
\bibitem{Tikhonov1}
Tikhonov A.N., Arsenin V.I. {\it Methods of regularizations of ill posed problem.}
Moscow, Nauka, (1975), in Russian.
\bibitem{Tsybakov1}
Tsybakov A.B. {\it  Introduction  a l estimation nonparametrique. } Springer;
(2004), New York, London, MR 2013911, in French.


\end{thebibliography}
\end{document}